\documentclass[amstex,12pt,russian,amssymb]{article}

\usepackage{mathtext}
\usepackage[cp1251]{inputenc}
\usepackage[T2A]{fontenc}
\usepackage[russian]{babel}
\usepackage[dvips]{graphicx}
\usepackage{amsmath}
\usepackage{amssymb}
\usepackage{amsxtra}
\usepackage{latexsym}
\usepackage{ifthen}

\begin{document}

\newcounter{lemma}
\newcommand{\lemma}{\par \refstepcounter{lemma}%
{\bf Лемма \arabic{lemma}.}}

\newcounter{corollary}
\newcommand{\corollary}{\par \refstepcounter{corollary}%
{\bf Следствие \arabic{corollary}.}}

\newcounter{remark}
\newcommand{\remark}{\par \refstepcounter{remark}%
{\bf Замечание \arabic{remark}.}}

\newcounter{theorem}
\newcommand{\theorem}{\par \refstepcounter{theorem}%
{\bf Теорема \arabic{theorem}.}}

\newcounter{example}
\newcommand{\example}{\par \refstepcounter{example}%
{\bf Пример \arabic{example}.}}

\newcounter{proposition}
\newcommand{\proposition}{\par \refstepcounter{proposition}%
{\bf Предложение \arabic{proposition}.}}

\renewcommand{\refname}{\centerline{\bf Список литературы}}

\newcommand{\proof}{{\it Доказательство.\,\,}}

\noindent УДК 517.5

{\bf Е.А.~Севостьянов} (Житомирский государственный университет
имени Ивана Франко, Институт прикладной математики и механики НАН Украины, г.~Славянск)

{\bf С.А.~Скворцов, А.П.~Довгопятый} (Житомирский государственный
университет имени Ивана Франко)

\medskip\medskip
{\bf Є.О.~Севостьянов} (Житомирський державний університет імені
Івана Фран\-ка, Інститут прикладної математики і мехініки НАН України, м.~Слов'янськ)

{\bf С.О.~Скворцов, О.П.~Довгопятий} (Житомирський державний
університет імені Івана Фран\-ка)

\medskip\medskip
{\bf E.A.~Sevost'yanov} (Zhytomyr Ivan Franko State University,
Institute of Applied Ma\-the\-ma\-tics and Mechanics of NAS of
Ukraine, Slov'yans'k)

{\bf S.O.~Skvortsov, O.P.~Dovhopiatyi} (Zhytomyr Ivan Franko State
University)

\medskip
{\bf Об отображениях, удовлетворяющих обратному неравенству Полецкого}

{\bf Про відображення, що задовольняють обернену нерівність Полецького}

{\bf On mappings satisfying the inverse Poletsky inequality}

\medskip\medskip
Изучено локальное и граничное поведение отображений,
удовлетворяющих одной оценке искажения модуля семейств кривых.
В частности, получены условия, при которых семейства
указанных отображений равностепенно непрерывны внутри и на границе области.

\medskip\medskip
Досліджено локальну і межову поведінку відображень,
які задовольняють одну оцінку спотворення модуля сімей кривих.
Зокрема, отримано умови, за яких сім'ї вказаних відоражень є одностайно неперервними всередині і на межі області.

\medskip\medskip
We have studied the local and boundary behavior of mappings satisfying one estimate of the distortion of the modulus of families of paths.
In particular, we have obtained conditions under which the families of the indicated mappings are equicontinuous inside and on the boundary of a domain.

\newpage
{\bf 1. Введение.} Хорошо известно, что отображения с ограниченным искажением (квазирегулярные отображения)
удовлетворяют в своей области определения соотношению вида
\begin{equation}\label{eq2}
M(\Gamma)\leqslant N(f, A)\cdot K\cdot M(f(\Gamma))\,,
\end{equation}
где $M$ -- модуль семейств кривых $\Gamma$ в области $D,$ $$N(y, f, A)\,=\,{\rm card}\,\left\{x\in A:
f(x)=y\right\}\,, \qquad N(f, A)\,=\,\sup\limits_{y\in{\Bbb
R}^n}\,N(y, f, A)\,,$$
$A$ -- произвольное борелевское множество в $D,$ а $K\geqslant 1$ --
некоторая постоянная, которая может быть вычислена как $K={\rm ess
\sup}\, K_O(x, f),$ $K_O(x, f)=\Vert f^{\,\prime}(x)\Vert^n/J(x, f)$
при $J(x, f)\ne 0;$ $K_O(x, f)=1$ при $f^{\,\prime}(x)=0,$ и $K_O(x,
f)=\infty$ при $f^{\,\prime}(x)\ne 0,$ но $J(x, f)=0$ (см., напр.,
\cite[теорема~3.2]{MRV$_1$} либо \cite[теорема~6.7, гл.~II]{Ri}).
Нечто аналогичное имеет место и для отображений с неограниченной
характеристикой. В частности, для так называемых отображений с
конечным искажением длины установлены оценки вида
\begin{equation}\label{eq1C}
M(\Gamma)\le \int\limits_{f(E)}K_I\left(y, f^{\,-1},
E\right)\cdot\rho_*^n(y)\,dm(y)\,,
\end{equation}
где $E$ -- произвольное измеримое подмножество области $D,$ $\Gamma$
-- произвольное семейство кривых в $E$ и $\rho_*(y)\in {\rm
adm\,}f(\Gamma)$ (см., напр., \cite[теорема~8.5]{MRSY}). В настоящей
заметке основным объектом изучения являются отображения,
удовлетворяющие некоторому более общему неравенству по сравнению
с~(\ref{eq1C}). Обратимся к определениям. Пусть $y_0\in {\Bbb R}^n,$
$0<r_1<r_2<\infty$ и
\begin{equation}\label{eq1**}
A(y_0, r_1,r_2)=\left\{ y\,\in\,{\Bbb R}^n:
r_1<|y-y_0|<r_2\right\}\,.\end{equation}
Для заданных множеств $E,$ $F\subset\overline{{\Bbb R}^n}$ и области
$D\subset {\Bbb R}^n$ обозначим через $\Gamma(E,F,D)$ семейство всех
кривых $\gamma:[a,b]\rightarrow \overline{{\Bbb R}^n}$ таких, что
$\gamma(a)\in E,\gamma(b)\in\,F$ и $\gamma(t)\in D$ при $t \in (a,
b).$ Если $f:D\rightarrow {\Bbb R}^n$ -- заданное отображение,
$y_0\in f(D)$ и $0<r_1<r_2<\infty,$ то через $\Gamma_f(y_0, r_1,
r_2)$ мы обозначим семейство всех кривых $\gamma$ в области $D$
таких, что $f(\gamma)\in \Gamma(S(y_0, r_1), S(y_0, r_2),
A(y_0,r_1,r_2)).$ Пусть $Q:{\Bbb R}^n\rightarrow [0, \infty]$ --
измеримая по Лебегу функция.  Будем говорить, что {\it $f$
удовлетворяет обратному неравенству Полецкого} в точке $y_0\in
f(D),$ если соотношение
\begin{equation}\label{eq2*A}
M(\Gamma_f(y_0, r_1, r_2))\leqslant \int\limits_{f(D)} Q(y)\cdot
\eta^n (|y-y_0|)\, dm(y)
\end{equation}
выполняется для произвольной измеримой по Лебегу функции $\eta:
(r_1,r_2)\rightarrow [0,\infty ]$ такой, что
\begin{equation}\label{eqA2}
\int\limits_{r_1}^{r_2}\eta(r)\, dr\geqslant 1\,.
\end{equation}
По поводу сравнения~(\ref{eq2*A}) с классическим неравенством
Полецкого укажем на~\cite[теорема~1]{Pol}. Отметим, что первым
автором установлены открытость и дискретность отображений
вида~(\ref{eq2*A}) при определённых условиях на функцию $Q,$ см.,
напр., \cite{Sev$_3$}. В более общем случае, относящемуся к данной
статье, выполнение этих свойств не гарантировано. Отметим также, что
равностепенная непрерывность гомеоморфизмов с условием~(\ref{eq2*A})
при несколько менее общих ограничениях на области и соответствующие
отображения подробно изучена в~\cite{SevSkv$_3$}. Данная работа
преимущественно относится к отображениям с ветвлением.

\medskip
Для области $D\subset {\Bbb R}^n,$ $n\geqslant 2,$ и измеримой по
Лебегу функции $Q:{\Bbb R}^n\rightarrow [0, \infty]$ обозначим через
$\frak{F}_Q(D)$ семейство всех открытых дискретных отображений
$f:D\rightarrow {\Bbb R}^n$ таких, что соотношение~(\ref{eq2*A})
выполнено для каждой точки $y_0\in f(D).$ Основные результаты,
полученные в статье, следующие.

\medskip
\begin{theorem}\label{th1}
{\sl\, Пусть $Q\in L^1({\Bbb R}^n).$ Тогда найдётся постоянная
$C_n>0,$ зависящая только от размерности пространства $n,$ такая что
для каждого $x_0\in D$ и каждого $0<2r_0<{\rm dist}\,(x_0,
\partial D)$ имеет место неравенство
\begin{equation}\label{eq2C}
|f(x)-f(x_0)|\leqslant\frac{C_n\cdot (\Vert
Q\Vert_1)^{1/n}}{\log^{1/n}\left(1+\frac{r_0}{|x-x_0|}\right)}\quad\forall\,\,x\in
B(x_0, r_0)\,,\quad \forall\,\, f\in \frak{F}_Q(D)\,,
\end{equation}
где $\Vert Q\Vert_1$ -- норма функции $Q$ в $L^1({\Bbb R}^n).$ В
частности, семейство $\frak{F}_Q(D)$ равностепенно непрерывно в
$D.$}
\end{theorem}

Напомним, что область $D\subset {\Bbb R}^n$ называется {\it локально
связной в точ\-ке} $x_0\in\partial D,$ если для любой окрестности
$U$ точки $x_0$ найдется окрестность $V\subset U$ точки $x_0$ такая,
что $V\cap D$ связно. Область $D$ локально связна на $\partial D,$
если $D$ локально связна в каждой точке $x_0\in\partial D.$ Граница
области $D$ называется {\it слабо плоской} в точке $x_0\in
\partial D,$ если для каждого $P>0$ и для любой окрестности $U$
точки $x_0$ найдётся окрестность $V\subset U$ этой же точки такая,
что $M(\Gamma(E, F, D))>P$ для произвольных континуумов $E, F\subset
D,$ пересекающих $\partial U$ и $\partial V.$ Граница области $D$
называется слабо плоской, если соответствующее свойство выполнено в
каждой точке границы $D.$

\medskip
Для числа $\delta>0,$ областей $D, D^{\,\prime}\subset {\Bbb R}^n,$
$n\geqslant 2,$ континуума $A\subset D^{\,\prime}$ и произвольной
измеримой по Лебегу функции $Q:D^{\,\prime}\rightarrow [0, \infty]$
обозначим через ${\frak S}_{\delta, A, Q }(D, D^{\,\prime})$
семейство всех открытых дискретных и замкнутых отображений $f$
области $D$ на $D^{\,\prime},$ удовлетворяющих условию~(\ref{eq2*A})
для каждого $y_0\in f(D),$ при этом,  $h(f^{\,-1}(A),
\partial D)=\inf\limits_{x\in f^{\,-1}(A), y\in \partial D}h(x, y)\geqslant~\delta.$ Справедливо
следующее утверждение.

\medskip
\begin{theorem}\label{th2}
{\sl Предположим, что область $D$ имеет слабо плоскую границу. Если
$Q\in L^1(D^{\,\prime}),$ и область $D^{\,\prime}$ локально связна
на своей границе, то каждое $f\in{\frak S}_{\delta, A, Q }(D,
D^{\,\prime})$ продолжается по непрерывности до отображения
$\overline{f}:\overline{D}\rightarrow \overline{D^{\,\prime}},$
$\overline{f}|_D=f,$ при этом,
$\overline{f}(\overline{D})=\overline{D^{\,\prime}}$ и семейство
${\frak S}_{\delta, A, Q }(\overline{D}, \overline{D^{\,\prime}}),$
состоящее из всех продолженных отображений
$\overline{f}:\overline{D}\rightarrow \overline{D^{\,\prime}},$
равностепенно непрерывно в $\overline{D}.$ }
\end{theorem}

\medskip
\begin{remark}\label{rem1}
В теореме~\ref{th1} равностепенную непрерывность следует понимать
относительно евклидовой метрики, а именно, для каждого
$\varepsilon>0$ найдётся $\delta=\delta(\varepsilon, x_0)>0$ такое,
что условие $|x-x_0|<\delta$ влечёт неравенство
$|f(x)-f(x_0)|<\varepsilon$ для всех $f\in\frak{F}_Q(D).$ В
теореме~\ref{th2} равностепенную непрерывность следует понимать
относительно хордальной метрики, т.е., для каждого $\varepsilon>0$
найдётся $\delta=\delta(\varepsilon, x_0)>0$ такое, что условие
$h(x, x_0)<\delta$ влечёт неравенство $h(\overline{f}(x),
\overline{f}(x_0))<\varepsilon$ при всех $\overline{f}\in{\frak
S}_{\delta, A, Q }(\overline{D}, \overline{D^{\,\prime}}).$
\end{remark}

\medskip
{\bf 2. О равностепенной непрерывности отображений во внутренних
точках.}

\begin{figure}[h]
\centerline{\includegraphics[scale=0.7]{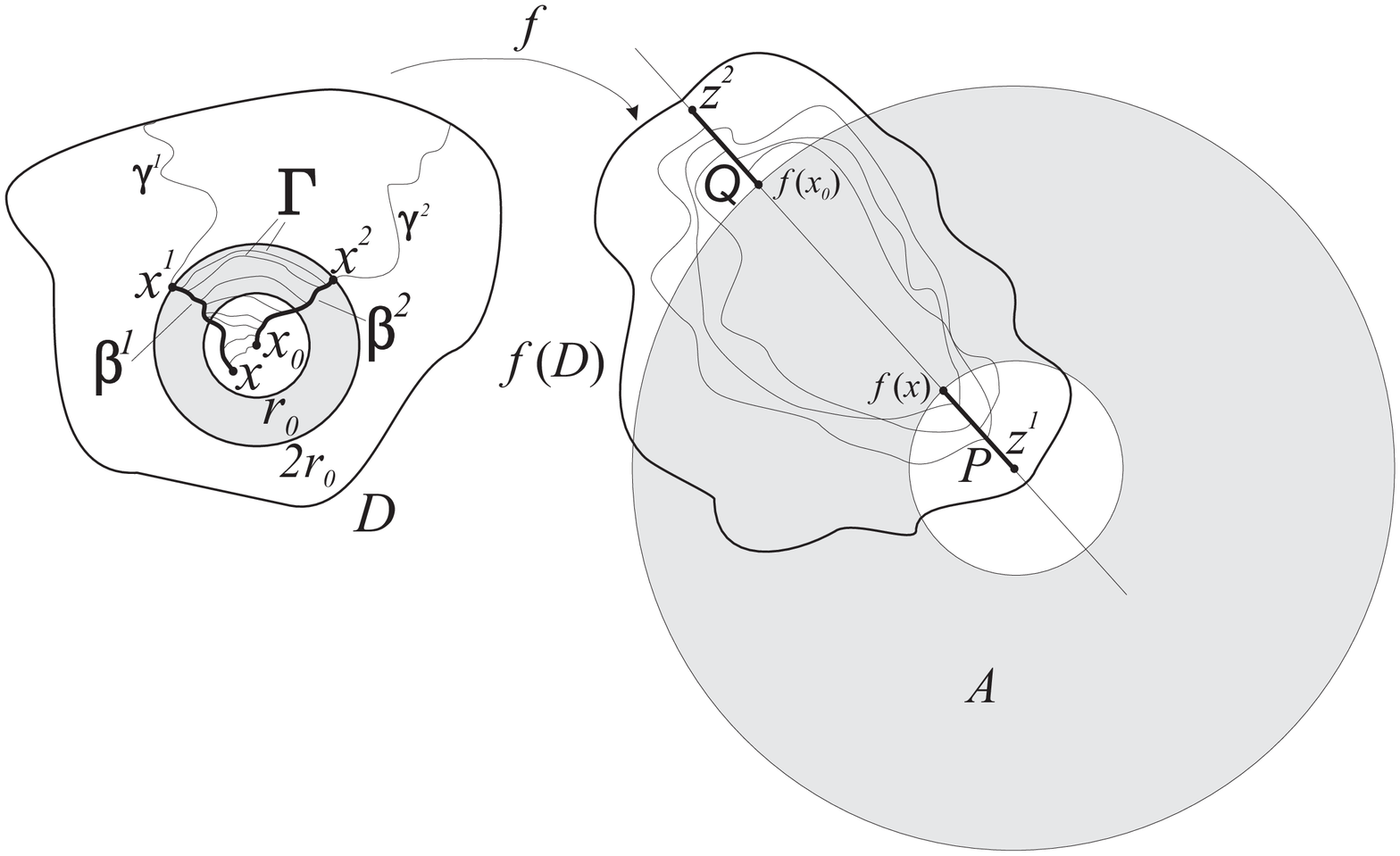}} \caption{К
доказательству теоремы~\ref{th1}}\label{fig2}
\end{figure}

\medskip
{\it Доказательство теоремы~\ref{th1}.} Зафиксируем $x_0\in D,$
$0<r_0<{\rm dist}\,(x_0,
\partial D)$ и $f\in
\frak{F}_Q(D).$ Рассмотрим $x\in B(x_0, r)$ и положим
\begin{equation}\label{eq13***}
|f(x)-f(x_0)|:=\varepsilon_0\,.
\end{equation}
Проведём через точки $f(x)$ и $f(x_0)$ прямую
$r=r(t)=f(x_0)+(f(x)-f(x_0))t,$ $-\infty<t<\infty$
(см.~рисунок~\ref{fig2}).
Пусть $\gamma^1:[1, c)\rightarrow D,$ $1<c\leqslant \infty$ --
максимальное $f$-поднятие луча $r=r(t),$ $t\geqslant 1,$ с началом в
точке $x,$ существующее ввиду~\cite[лемма~3.12]{MRV$_3$}. По этой же
лемме
\begin{equation}\label{eq3}
h(\gamma^1(t), \partial D)\rightarrow 0
\end{equation}
при $t\rightarrow c-0,$ где $h(A, B)=\sup\limits_{x\in A,y\in B}h(x,
y)$ и $h(x, y)$ -- хордальное расстояние между точками $x, y\in
\overline{{\Bbb R}^n}$ (см.~\cite[определение~12.1]{Va}).

\medskip
Аналогично, пусть $\gamma^2:(d, 0]\rightarrow D,$ $-\infty\leqslant
d<0$ -- максимальное $f$-поднятие луча $r=r(t),$ $t\leqslant 0,$ с
началом в точке $x_0,$ существующее
ввиду~\cite[лемма~3.12]{MRV$_3$}. Также, как и в~(\ref{eq3}) мы
имеем, что
\begin{equation}\label{eq4}
h(\gamma^2(t), \partial D)\rightarrow 0
\end{equation}
при $t\rightarrow d-0.$ Из соотношений~(\ref{eq3}) и~(\ref{eq4}), с
учётом~\cite[теорема~1.I.5.46]{Ku}, вытекает существование чисел
$1\leqslant t^1<c$ и $1\leqslant t^2<d$ и элементов
$x^1:=\gamma^1(t^1)$ и $x^2:=\gamma^2(t^2)\in S(x_0, 2r_0).$ Не
ограничивая общности, можно считать, что $\gamma^1(t)\in B(x_0,
2r_0)$ при всех $t\in [1, t^1]$ и $\gamma^2(t)\in B(x_0, 2r_0)$ при
всех $t\in [t^2, 0].$
Пусть
$$\beta^1:=\gamma^1|_{[1, t^1]},\qquad \beta^2:=\gamma^2|_{[1,
t^2]},\qquad \Gamma:=(|\beta^1|, |\beta^2|, B(x_0, 2r_0))\,.$$
Тогда с одной стороны ввиду~\cite[лемма~4.3]{Vu$_2$}
\begin{equation}\label{eq7A}
M(\Gamma)\geqslant (1/2)\cdot M(\Gamma(|\beta^1|, |\beta^2|, {\Bbb
R}^n))\,,
\end{equation}
а с другой стороны, по~\cite[лемма~7.38]{Vu$_3$}
\begin{equation}\label{eq7B}
M(\Gamma(|\beta^1|, |\beta^2|, {\Bbb R}^n))\geqslant
c_n\cdot\log\left(1+\frac1m\right)\,,
\end{equation}
где $c_n>0$ -- некоторая постоянная, зависящая только от $n,$
$$m=\frac{d(|\beta^1|, |\beta^2|)}{\min\{d(|\beta^1|),
d(|\beta^2|)\}},\quad d(|\beta^1|, |\beta^2|)=\inf\limits_{x\in
|\beta^1|, y\in|\beta^2|}|x-y|\,.$$
Заметим, что ${\rm diam}\,(|\beta^i|)=\sup\limits_{x,
y\in|\beta^i|}|x-y|\geqslant r_0,$ $i=1,2.$ Тогда
объединяя~(\ref{eq7A}) и~(\ref{eq7B}) и учитывая, что $d(|\beta^1|,
|\beta^2|)\leqslant |x-x_0|,$  мы получаем, что
\begin{equation}\label{eq7C}
M(\Gamma)\geqslant \widetilde{c_n}\cdot
\log\left(1+\frac{r_0}{d(|\beta^1|, |\beta^2|)}\right)\geqslant
\widetilde{c_n}\cdot \log\left(1+\frac{r_0}{|x-x_0|}\right)\,,
\end{equation}
где $\widetilde{c_n}>0$ -- некоторая постоянная, зависящая только от
$n.$

\medskip
Установим теперь оценку сверху для $M(\Gamma).$ Пусть $P$ -- часть
прямой $r(t),$ расположенная между точками $f(x)$ и
$z^1:=f(x^1)=f(\gamma^1(t^1)),$ а $Q$ -- часть прямой $r(t),$
расположенная между точками $f(x_0)$ и
$z^2:=f(x^2)=f(\gamma^2(t^2)).$ Положим
$$A:=A(z^1, \varepsilon^1, \varepsilon^2)=\{x\in {\Bbb R}^n: \varepsilon^1<|x-z^1|<\varepsilon^2\}\,,$$
где $\varepsilon^1:=|f(x)- z^1|,$ $\varepsilon^2:=|f(x_0)-z^1|.$
Покажем, что
\begin{equation}\label{eq1A}
f(\Gamma)>\Gamma(S(z^1, \varepsilon^1), S(z^1, \varepsilon^2), A)\,.
\end{equation}
В самом деле, пусть $\gamma\in\Gamma.$ Тогда $f(\gamma)\in
f(\Gamma),$ $f(\gamma)=f(\gamma(s)):[0, 1]\rightarrow {\Bbb R}^n,$
$f(\gamma(0))\in P,$ $f(\gamma(1))\in Q$ и $f(\gamma(s))\in f(D)$
при $0<s<1.$ Пусть $q>1$ -- число, такое что
$$z^1=f(x_0)+(f(x)-f(x_0))q\,.$$ Поскольку
$f(\gamma(0))\in P,$ найдётся $1\leqslant t\leqslant q$ такое, что
$f(\gamma(0))=f(x_0)+(f(x)-f(x_0))t.$ Следовательно,
$$|f(\gamma(0))-z^1|=|(f(x)-f(x_0))(q-t)|\leqslant $$
\begin{equation}\label{eq2A}\leqslant
|(f(x)-f(x_0))(q-1)|=|(f(x)-f(x_0))q+f(x_0)-f(x))|=\end{equation}
$$=|f(x)-z^1|=\varepsilon^1\,.$$

С другой стороны, поскольку $f(\gamma(1))\in Q,$ найдётся
$p\leqslant 0$ такое, что
$$f(\gamma(1))=f(x_0)+(f(x)-f(x_0))p\,.$$ В этом случае,
мы получим, что
$$
|f(\gamma(1))-z^1|=|(f(x)-f(x_0))(q-p)|\geqslant $$
\begin{equation}\label{eq3A}
\geqslant |(f(x)-f(x_0))q|=|(f(x)-f(x_0))q +f(x_0)-f(x_0)|=
\end{equation}
$$=|f(x_0)-z^1|=\varepsilon^2\,.$$
Заметим, что
\begin{equation}\label{eq5B}
|f(x_0)- f(x)|+\varepsilon^1=|f(x_0)- f(x)|+|f(x)-z^1|=
|z^1-f(x_0)|=\varepsilon^2\,,
\end{equation}
и, значит, $\varepsilon^1<\varepsilon^2.$ Тогда из~(\ref{eq3A})
вытекает, что
\begin{equation}\label{eq4A}
|f(\gamma(1))-z^1|>\varepsilon^1\,.
\end{equation}
Пусть $f(\gamma(0))\not\in S(z^1, \varepsilon^1),$ тогда
из~(\ref{eq2A}) и (\ref{eq4A}) вытекает, что $|f(\gamma)|\cap B(z^1,
\varepsilon^1)\ne\varnothing\ne (f(D)\setminus B(z^1,
\varepsilon^1))\cap|f(\gamma)|.$ В таком случае, в
силу~\cite[теорема~1.I.5.46]{Ku} найдётся $t_1\in (0, 1)$ такое, что
$f(\gamma(t_1))\in S(z^1, \varepsilon^1).$ Без ограничения общности,
мы можем считать, что $f(\gamma(t))\not\in B(z^1, \varepsilon^1)$
при $t\in (t_1, 1).$ Положим $\alpha^1:=f(\gamma)|_{[t_1, 1]}.$

\medskip
С другой стороны, поскольку $\varepsilon^1<\varepsilon^2$ и
$f(\gamma)_1(t_1)\in S(z^1, \varepsilon^1),$ мы получим, что
$|\alpha^1|\cap B(z^1, \varepsilon^2)\ne\varnothing.$ По
соотношению~(\ref{eq3A}) мы получим, что $(f(D)\setminus B(z^1,
\varepsilon^2))\cap |\alpha^1|\ne\varnothing.$ Следовательно,
по~\cite[теорема~1.I.5.46]{Ku} существует $t_2\in [t_1, 1)$ такое,
что $\alpha^1(t_2)\in S(z^1, \varepsilon^2).$

Без ограничения общности, мы можем считать, что $f(\gamma)_1(t)\in
B(z^1, \varepsilon^2)$ при $t\in (t_1, t_2).$ Положим
$\alpha^2:=\alpha^1|_{[t_1, t_2]}.$ Тогда $f(\gamma)>\alpha^2$ и
$\alpha^2\in \Gamma(S(z^1, \varepsilon^1), S(z^1, \varepsilon^2),
A).$ Таким образом, соотношение~(\ref{eq1A}) доказано.

\medskip
Из~(\ref{eq1A}) вытекает, что $\Gamma>\Gamma_{f}(z^1, \varepsilon^1,
\varepsilon^2).$ Теперь полагаем
$$\eta(t)= \left\{
\begin{array}{rr}
\frac{1}{\varepsilon_0}, & t\in [\varepsilon^1, \varepsilon^2],\\
0,  &  t\not\in [\varepsilon^1, \varepsilon^2]\,.
\end{array}
\right. $$
Заметим, что $\eta$ удовлетворяет соотношению~(\ref{eqA2}) при
$r_1=\varepsilon^1$ и $r_2=\varepsilon^2.$ В самом деле,
из~(\ref{eq13***}) и (\ref{eq5B}) вытекает, что
$$r_1-r_2=\varepsilon^2-\varepsilon^1=|f(x_0)-z^1|-|f(x)-
z^1|=$$$$=|f(x)-f(x_0)|=\varepsilon_0\,.$$ Тогда
$\int\limits_{\varepsilon^1}^{\varepsilon^2}\eta(t)\,dt=(1/\varepsilon_0)\cdot
(\varepsilon^2-\varepsilon^1)\geqslant 1.$ По
неравенству~(\ref{eq1A}), рассматриваемому в точке $z^1$ и
положенного в основу определения семейства $\frak{F}_Q(D)$ мы
получим, что
$$M(\Gamma)\leqslant M(\Gamma_{f}(z^1, \varepsilon^1, \varepsilon^2)) \leqslant$$
\begin{equation}\label{eq14***}
\leqslant \frac{1}{\varepsilon_0^n}\int\limits_{{\Bbb R}^n}
Q(y)\,dm(y)=\frac{\Vert Q\Vert_1}{{|f(x)-f(x_0)|}^n}\,.
\end{equation}
Из~(\ref{eq7C}) и (\ref{eq14***}) вытекает, что
$$\widetilde{c_n}\cdot \log\left(1+\frac{r_0}{|x-x_0|}\right)\leqslant
\frac{\Vert Q\Vert_1}{{|f(x)-f(x_0)|}^n}\,.$$
Из последнего соотношения вытекает желанное
неравенство~(\ref{eq2C}), где
$C_n:={\widetilde{c_n}}^{-1/n}.$~$\Box$

\medskip
{\bf 3. Граничное поведение отображений с обратным неравенством
Полецкого.} Отметим некоторые известные утверждения о продолжении
гомеоморфизмов с условием~(\ref{eq2*A}) на границу области, см.,
напр., \cite[лемма~5.20, следствие~5.23]{MRSY$_1$}, \cite[лемма~6.1,
теорема~6.1]{RS} и \cite[лемма~5, теорема~3]{Sm}. Наша ближайшая
цель -- получить аналогичный результат для отображений, допускающих
ветвление. Имеет место следующее утверждение.

\medskip
\begin{theorem}\label{th3}
{\it Пусть $D\subset {\Bbb R}^n,$ $n\geqslant 2,$ -- область,
имеющая слабо плоскую границу, а область $D^{\,\prime}\subset {\Bbb
R}^n$ локально связна на своей границе. Предположим, $f$ -- открытое
дискретное и замкнутое отображение области $D$ на $D^{\,\prime},$
удовлетворящее соотношению~(\ref{eq2*A}) в каждой точке $y_0\in
f(D),$ где $Q\in L^1(D^{\,\prime}).$ Тогда отображение $f$
продолжается по непрерывности до отображения
$\overline{f}:\overline{D}\rightarrow\overline{D^{\,\prime}},$ при
этом, $\overline{f}(\overline{D})=\overline{D^{\,\prime}}.$}
\end{theorem}

\medskip
\begin{proof}
Зафиксируем произвольным образом точку $x_0\in\partial D.$
Необходимо показать возможность непрерывного продолжения отображения
$f$ в точку $x_0.$ Применяя при необходимости мёбиусово
преобразование $\varphi:\infty\mapsto 0$ и учитывая инвариантность
модуля $M$ в левой части соотношения~(\ref{eq2*A})
(см.~\cite[теорема~8.1]{Va}), мы можем считать, что $x_0\ne\infty.$

\medskip
Предположим, что заключение о непрерывном продолжении отображения
$f$ в точку $x_0$ не является верным. Тогда найдётся не менее двух
последовательностей $x_i, y_i\in D,$ $i=1,2,\ldots ,$ таких, что
$x_i, y_i\rightarrow x_0$ при $i\rightarrow\infty,$ при этом,
$h(f(x_i), f(y_i))\geqslant a>0$ при некотором $a>0$ и всех $i\in
{\Bbb N},$ где $h$ -- хордальная метрика,
см.~\cite[определение~12.1]{Va}. В силу компактности пространства
$\overline{{\Bbb R}^n}$ мы можем считать, что последовательности
$f(x_i)$ и $f(y_i)$ сходятся при $i\rightarrow\infty$ к $z_1$ и
$z_2,$ соответственно, причём $z_1\ne\infty.$ Так как отображение
$f$ замкнуто, то оно сохраняет границу области,
см.~\cite[теорема~3.3]{Vu$_1$}, поэтому $z_1, z_2\in \partial
D^{\,\prime}.$ Поскольку область $D^{\,\prime}$ локально связна на
своей границе, существуют непересекающиеся окрестности $U_1$ и $U_2$
точек $z_1$ и $z_2$ такие, что $W_1=D^{\,\prime}\cap U_1$ и
$W_2=D^{\,\prime}\cap U_2$ являются связными. Можно считать, что
$W_1$ и $W_2$ линейно связные, поскольку в качестве окрестностей
$U_1$ и $U_2$ можно выбрать открытые множества (см.,
напр.,~\cite[предложение~13.2]{MRSY}). См. рисунок~\ref{fig1} по
поводу иллюстрации.
\begin{figure}[h]
\centerline{\includegraphics[scale=0.65]{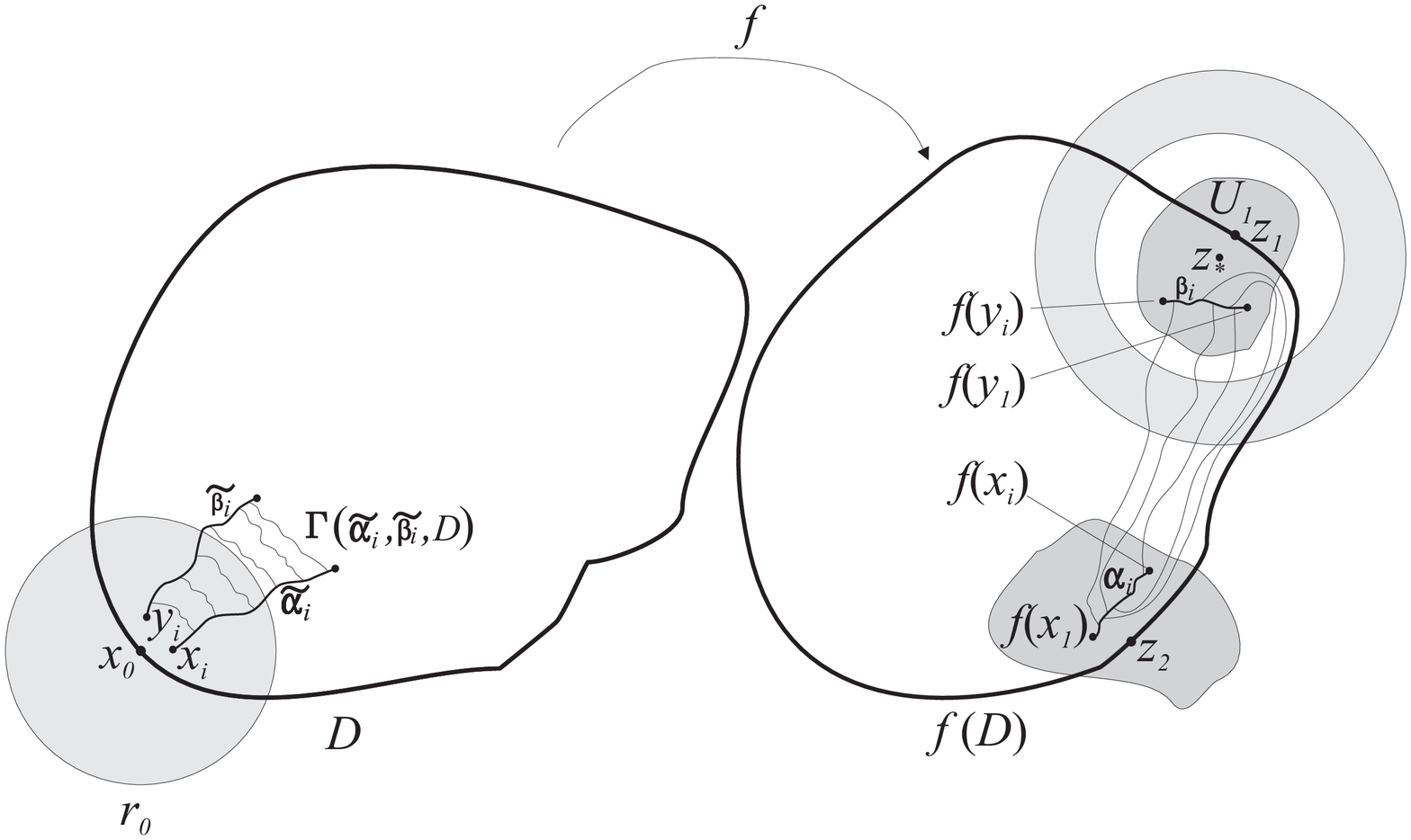}} \caption{К
доказательству теоремы~\ref{th3}}\label{fig1}
\end{figure}
Можно считать, что
\begin{equation}\label{eq8}
U_1\subset B(z_*, R_0),\qquad \overline{B(z_*, 2R_0)}\cap
\overline{U_2}=\varnothing\,,\qquad R_0>0\,,
\end{equation}
где $z_*\in D^{\,\prime}$ -- некоторая точка, достаточно близкая к
$z_1.$ Мы также можем считать, что $f(x_i)\in W_1$ и $f(y_i)\in W_2$
при всех $i=1,2,\ldots .$ Соединим точки $f(x_i)$ и $f(x_1)$ кривой
$\alpha_i:[0, 1]\rightarrow D^{\,\prime},$ а точки $f(y_i)$ и
$f(y_1)$ кривой $\beta_i:[0, 1]\rightarrow D^{\,\prime}$ так, что
$|\alpha_i|\subset W_1$ и $|\beta_i|\subset W_2$ при $i=1,2,\ldots
.$ Пусть $\widetilde{\alpha_i}:[0, 1]\rightarrow D^{\,\prime}$ и
$\widetilde{\beta_i}:[0, 1]\rightarrow D^{\,\prime}$ -- полные
поднятия кривых $\alpha_i$ и $\beta_i$ с началами в точках $x_i$ и
$y_i,$ соответственно (эти поднятия существуют
ввиду~\cite[лемма~3.7]{Vu$_1$}). Заметим, что у точек $f(x_1)$ и
$f(y_1)$ в области $D$ может быть не более конечного числа
прообразов при отображении $f,$ см.~\cite[лемма~3.2]{Vu$_1$}. Тогда
найдётся $r_0>0$ такое, что $\widetilde{\alpha_i}(1),
\widetilde{\beta_i}(1)\in D\setminus B(x_0, r_0)$ при всех
$i=1,2,\ldots .$ Поскольку граница области $D$ является слабо
плоской, для каждого $P>0$ найдётся $i=i_P\geqslant 1$ такое, что
\begin{equation}\label{eq7}
M(\Gamma(|\widetilde{\alpha_i}|, |\widetilde{\beta_i}|,
D))>P\qquad\forall\,\,i\geqslant i_P\,.
\end{equation}
Покажем, что условие~(\ref{eq7}) противоречит определению
отображения $f$ в~(\ref{eq2*A}). В самом деле, в силу
соотношений~(\ref{eq8}) и ввиду~\cite[теорема~1.I.5.46]{Ku}
\begin{equation}\label{eq9}
f(\Gamma(|\widetilde{\alpha_i}|, |\widetilde{\beta_i}|,
D))>\Gamma(S(z_*, R_0), S(z_*, 2R_0), A(z_*, R_0, 2R_0))\,.
\end{equation}
Из~(\ref{eq9}) вытекает, что
\begin{equation}\label{eq10}
\Gamma(|\widetilde{\alpha_i}|, |\widetilde{\beta_i}|, D)
>\Gamma_f(z_0, R_0, 2R_0)\,.
\end{equation}
Положим $\eta(t)= \left\{
\begin{array}{rr}
\frac{1}{R_0}, & t\in [R_0, 2R_0],\\
0,  &  t\not\in [R_0, 2R_0]
\end{array}
\right. .$
Заметим, что $\eta$ удовлетворяет соотношению~(\ref{eqA2}) при
$r_1=R_0$ и $r_2=2R_0.$ Тогда из~(\ref{eq10}) и (\ref{eq2*A}) мы
получим, что
$$M(\Gamma(|\widetilde{\alpha_i}|, |\widetilde{\beta_i}|, D))
\leqslant M(\Gamma_f(z_0, R_0, 2R_0)) \leqslant$$
\begin{equation}\label{eq11}
\leqslant \frac{1}{R_0^n}\int\limits_{D^{\,\prime}}
Q(y)\,dm(y):=c<\infty\,,
\end{equation}
поскольку~$Q\in L^1(D).$ Однако, соотношение~(\ref{eq11})
противоречит условию~(\ref{eq7}). Полученное противоречие
опровергает предположение об отсутствии предела у отображения $f$ в
точке $x_0.$

Осталось проверить равенство
$\overline{f}(\overline{D})=\overline{D^{\,\prime}}.$ Очевидно,
$\overline{f}(\overline{D})\subset\overline{D^{\,\prime}}.$ Покажем,
что $\overline{D^{\,\prime}}\subset \overline{f}(\overline{D}).$ В
самом деле, пусть $y_0\in \overline{D^{\,\prime}},$ тогда либо
$y_0\in D^{\,\prime},$ либо $y_0\in \partial D^{\,\prime}.$ Если
$y_0\in D^{\,\prime},$ то $y_0=f(x_0)$ и $y_0\in
\overline{f}(\overline{D}),$ поскольку по условию $f$ -- отображение
области $D$ на $D^{\,\prime}.$ Наконец, пусть $y_0\in \partial
D^{\,\prime},$ тогда найдётся последовательность $y_k\in
D^{\,\prime}$ такая, что $y_k=f(x_k)\rightarrow y_0$ при
$k\rightarrow\infty$ и $x_k\in D.$ В силу компактности пространства
$\overline{{\Bbb R}^n}$ мы можем считать, что $x_k\rightarrow x_0,$
где $x_0\in\overline{D}.$ Заметим, что $x_0\in \partial D,$ так как
отображение $f$ является открытым. Тогда $f(x_0)=y_0\in
\overline{f}(\partial D)\subset \overline{f}(\overline{D}).$ Теорема
полностью доказана.~$\Box$
\end{proof}

\medskip
{\bf 4. О равностепенной непрерывности семейств отображений в
замыкании области.}

\medskip
{\it Доказательство теоремы~\ref{th2}.} Пусть $f\in {\frak
S}_{\delta, A, Q }(D, D^{\,\prime}).$ По теореме~\ref{th3}
отображение $f$ продолжается до непрерывного отображения
$\overline{f}:\overline{D}\rightarrow \overline{D^{\,\prime}},$ при
этом, $\overline{f}(\overline{D})=\overline{D^{\,\prime}}.$
Равностепенная непрерывность семейства ${\frak S}_{\delta, A, Q
}(\overline{D}, \overline{D^{\,\prime}})$ в $D$ является
утверждением теоремы~\ref{th1}. Осталось установить равностепенную
непрерывность указанного семейства на $\partial D^{\,\prime}.$

\begin{figure}[h]
\centerline{\includegraphics[scale=0.7]{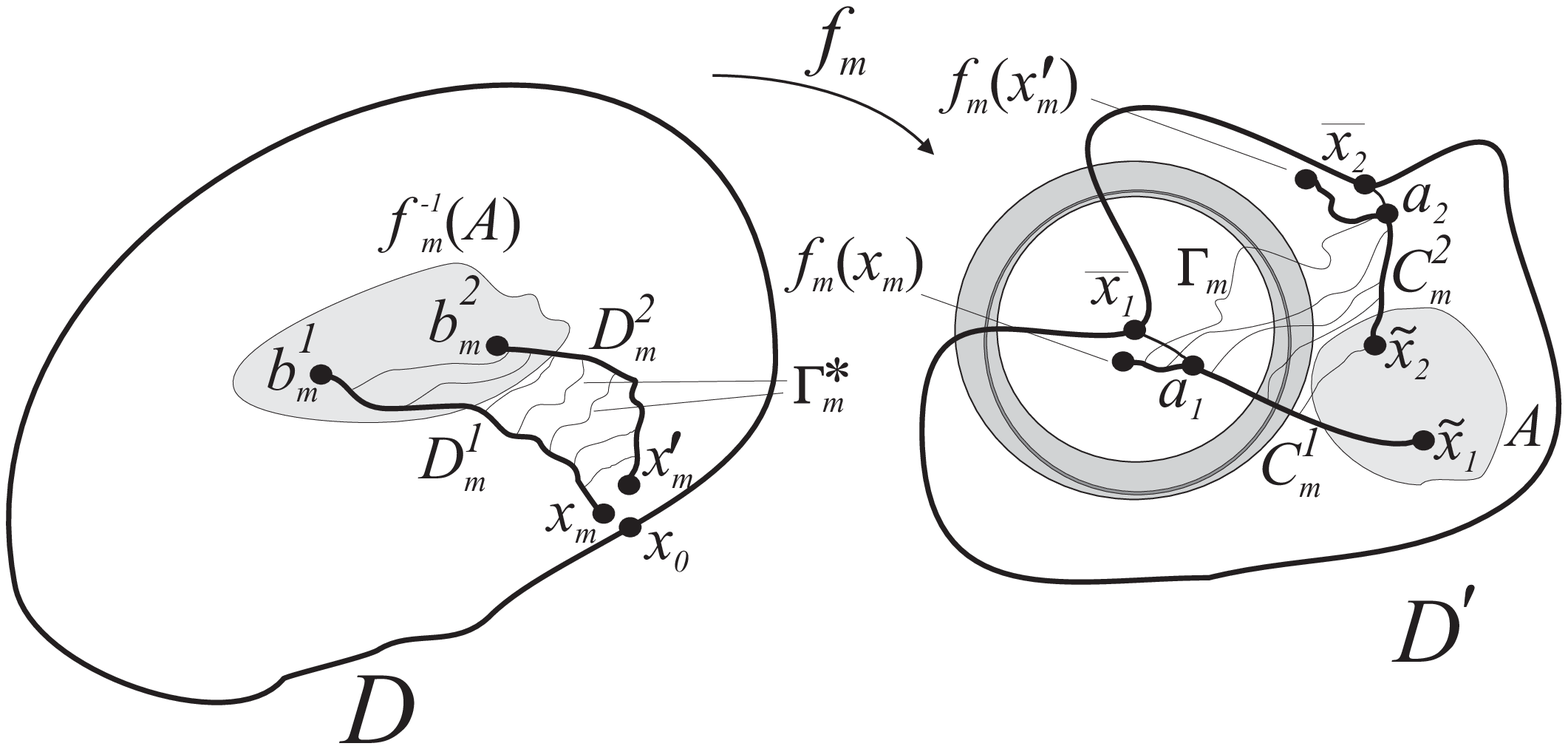}} \caption{К
доказательству теоремы~\ref{th2}}\label{fig4}
\end{figure}

Проведём доказательство от противного. Предположим, что найдётся
$x_0\in
\partial D,$ число $\varepsilon_0>0,$ последовательность $x_m\in \overline{D},$
сходящаяся к точке $x_0$ и соответствующие отображения
$\overline{f}_m\in {\frak S}_{\delta, A, Q }(\overline{D},
\overline{D})$ такие, что
\begin{equation}\label{eq12}
h(\overline{f}_m(x_m),\overline{f}_m(x_0))\geqslant\varepsilon_0,\quad
m=1,2,\ldots .
\end{equation}
Положим $f_m:=\overline{f}_m|_{D}.$ Поскольку $f_m$ имеет
непрерывное продолжение на $\partial D,$ можно считать, что $x_m\in
D.$ Следовательно, $\overline{f}_m(x_m)=f_m(x_m).$ Кроме того,
найдётся последовательность $x^{\,\prime}_m\in D$ такая, что
$x^{\,\prime}_m\rightarrow x_0$ при $m\rightarrow\infty$ и
$h(f_m(x^{\,\prime}_m),\overline{f}_m(x_0))\rightarrow 0$ при
$m\rightarrow\infty.$ Поскольку пространство $\overline{{\Bbb R}^n}$
компактно, мы можем считать, что $f_m(x_m)$ и $\overline{f}_m(x_0)$
сходятся при $m\rightarrow\infty.$ Пусть $f_m(x_m)\rightarrow
\overline{x_1}$ и $\overline{f}_m(x_0)\rightarrow \overline{x_2}$
при $m\rightarrow\infty.$ По непрерывности метрики в~(\ref{eq12}),
$\overline{x_1}\ne \overline{x_2}.$ Без ограничения общности, можно
считать, что $\overline{x_1}\ne\infty.$ Поскольку отображения $f_m$
замкнуты, то они сохраняют границу (см.~\cite[теорема~3.3]{Vu$_1$}),
поэтому $\overline{x_2}\in\partial D.$ Пусть $\widetilde{x_1}$ и
$\widetilde{x_2}$ различные точки континуума $A,$ ни одна из которых
не совпадает с $\overline{x_1}.$ По~\cite[лемма~2.1]{SevSkv$_2$} две
пары точек $\widetilde{x_1},$ $\overline{x_1}$ и $\widetilde{x_2},$
$\overline{x_2}$ можно соединить кривыми $\gamma_1:[0, 1]\rightarrow
\overline{D}$ и $\gamma_2:[0, 1]\rightarrow \overline{D}$ такими,
что $|\gamma_1|\cap |\gamma_2|=\varnothing,$ $\gamma_1(t),
\gamma_2(t)\in D$ при $t\in (0, 1),$ $\gamma_1(0)=\widetilde{x_1},$
$\gamma_1(1)=\overline{x_1},$ $\gamma_2(0)=\widetilde{x_2}$ и
$\gamma_2(1)=\overline{x_2}.$ Поскольку $D^{\,\prime}$ локально
связна  на $\partial D^{\,\prime},$ найдутся окрестности $U_1$ и
$U_2$ точек $\overline{x_1}$ и $\overline{x_2},$ чьи замыкания не
пересекаются, при этом, множества $W_i:=D^{\,\prime}\cap U_i$
линейно связны. Без ограничения общности можно считать, что
$\overline{U_1}\subset B(\overline{x_1}, \delta_0)$ и
\begin{equation}\label{eq12C}
\overline{B(\overline{x_1},
\delta_0)}\cap|\gamma_2|=\varnothing=\overline{U_2}\cap|\gamma_1|\,,
\quad \overline{B(\overline{x_1}, \delta_0)}\cap
\overline{U_2}=\varnothing\,.
\end{equation}
Можно также считать, что $f_m(x_m)\in W_1$ и $f_m(x^{\,\prime}_m)\in
W_2$ при всех $m\in {\Bbb N}.$ Пусть $a_1$ и $a_2$ -- две различные
точки, принадлежащие $|\gamma_1|\cap W_1$ и $|\gamma_2|\cap W_2.$
Пусть $0<t_1, t_2<1$ таковы, что $\gamma_1(t_1)=a_1$ и
$\gamma_2(t_2)=a_2.$ Соединим точки $a_1$ и $f_m(x_m)$ кривой
$\alpha_m:[t_1, 1]\rightarrow W_1$ такой, что $\alpha_m(t_1)=a_1$ и
$\alpha_m(1)=f_m(x_m).$ Аналогично, соединим $a_2$ и
$f_m(x^{\,\prime}_m)$ кривой $\beta_m:[t_2, 1]\rightarrow W_2,$
$\beta_m(t_2)=a_2$ и $\beta_m(1)=f_m(x^{\,\prime}_m)$ (см.
рисунок~\ref{fig4}).
Положим
$$C^1_m(t)=\quad\left\{
\begin{array}{rr}
\gamma_1(t), & t\in [0, t_1],\\
\alpha_m(t), & t\in [t_1, 1]\end{array} \right.\,,\qquad
C^2_m(t)=\quad\left\{
\begin{array}{rr}
\gamma_2(t), & t\in [0, t_2],\\
\beta_m(t), & t\in [t_2, 1]\end{array} \right.\,.$$
Пусть $D^1_m$ и $D^2_m$ -- полные поднятия кривых $|C^1_m|$ и
$|C^2_m|$ с началами в точках $x_m$ и $x^{\,\prime}_m,$
соответственно. Такие поднятия существуют по~\cite{Vu$_1$}, в
частности, в силу условия $h(f_m^{\,-1}(A),
\partial D)=\inf\limits_{x\in f_m^{\,-1}(A), y\in \partial D}h(x,
y)\geqslant~\delta,$ участвующего в определении класса ${\frak
S}_{\delta, A, Q }(D, D^{\,\prime}),$ концы кривых $D^1_m$ и
$D^2_m,$ обозначаемые в дальнейшем $b_m^1$ и $b_m^2,$ отстоят от
границы $D$ на расстояние, не меньшее $\delta.$

Как обычно, обозначим через $|C^1_m|$ и $|C^2_m|$ носители кривых
$C^1_m$ и $C^2_m,$ соответственно. Полагая
$$l_0=\min\{{\rm dist}\,(|\gamma_1|,
|\gamma_2|), {\rm dist}\,(|\gamma_1|, U_2)\}\,,$$
рассмотрим покрытие $A_0:=\bigcup\limits_{x\in |\gamma_1|}B(x,
l_0/4)$ кривой $|\gamma_1|$ посредством шаров. Поскольку
$|\gamma_1|$ -- компактное множество, можно выбрать конечное число
индексов $1\leqslant N_0<\infty$ и соответствующие точки
$z_1,\ldots, z_{N_0}\in |\gamma_1|$ так, что $|\gamma_1|\subset
B_0:=\bigcup\limits_{i=1}^{N_0}B(z_i, l_0/4).$
В этом случае,
$$|C^1_m|\subset U_1\cup |\gamma_1|\subset
\overline{B(\overline{x_1}, \delta_0)}\cup \bigcup\limits_{i=1}^{N_0}B(z_i, l_0/4)\,.$$
Пусть $\Gamma_m$ -- семейство всех кривых, соединяющих $|C^1_m|$ и
$|C^2_m|$ в $D^{\,\prime}.$ Тогда мы имеем, что
\begin{equation}\label{eq10C}
\Gamma_m=\bigcup\limits_{i=0}^{N_0}\Gamma_{mi}\,,
\end{equation}
где $\Gamma_{mi}$ -- семейство всех кривых $\gamma:[0, 1]\rightarrow
D^{\,\prime}$ таких, что $\gamma(0)\in B(z_i, l_0/4)\cap |C^1_m|$ и
$\gamma(1)\in |C_2^m|$ при $1\leqslant i\leqslant N_0.$ Аналогично,
пусть $\Gamma_{m0}$ -- семейство кривых $\gamma:[0, 1]\rightarrow
D^{\,\prime}$ таких, что $\gamma(0)\in B(\overline{x_1},
\delta_0)\cap |C^1_m|$ и $\gamma(1)\in |C_2^m|.$ В
силу~(\ref{eq12C}) найдётся $\sigma_0>\delta_0>0$ такое, что
$$
\overline{B(\overline{x_1},
\sigma_0)}\cap|\gamma_2|=\varnothing=\overline{U_2}\cap|\gamma_1|\,,
\quad \overline{B(\overline{x_1}, \sigma_0)}\cap
\overline{U_2}=\varnothing\,.$$
Рассуждая также, как при доказательстве теоремы~\ref{th1} и
используя~\cite[теорема~1.I.5.46]{Ku}, мы можем показать, что
$$\Gamma_{m0}>\Gamma(S(\overline{x_1}, \delta_0), S(\overline{x_1}, \sigma_0), A(\overline{x_1}, \delta_0,
\sigma_0))\,,$$
\begin{equation}\label{eq11C}
\Gamma_{mi}>\Gamma(S(z_i, l_0/4), S(z_i, l_0/2), A(z_i, l_0/4,
l_0/2))\,.
\end{equation}
Можно подобрать $x_*\in D^{\,\prime},$ $\delta_*>0$ и $\sigma_*>0$
такие, что $A(x_*, \delta_*, \sigma_*)\subset A(\overline{x_1},
\delta_0, \sigma_0),$ поэтому
\begin{equation}\label{eq1}
\Gamma(S(\overline{x_1}, \delta_0), S(\overline{x_1}, \sigma_0),
A(\overline{x_1}, \delta_0, \sigma_0))>\Gamma(S(x_*, \delta_*),
S(x_*, \sigma_*), A(x_*, \delta_*, \sigma_*))\,.
\end{equation}
Положим
$$\eta(t)= \left\{
\begin{array}{rr}
4/l_0, & t\in [l_0/4, l_0/2],\\
0,  &  t\not\in [l_0/4, l_0/2]\,,
\end{array}
\right. \qquad \eta_0(t)= \left\{
\begin{array}{rr}
1/(\sigma_*-\delta_*), & t\in [\delta_*, \sigma_*],\\
0,  &  t\not\in [\delta_*, \sigma_*]\,.
\end{array}
\right.$$
Обозначим $\Gamma^{\,*}_m:=\Gamma(|D_m^1|, |D_m^2|, D).$ Заметим,
что $f(\Gamma^{\,*}_m)\subset\Gamma_m.$ Тогда в силу (\ref{eq10C}),
(\ref{eq11C}) и~(\ref{eq1})
\begin{equation}\label{eq6A}
\Gamma^{\,*}_m>\left(\bigcup\limits_{i=1}^{N_0}\Gamma_{f_m}(z_i,
l_0/4, l_0/2)\right)\cup \Gamma_{f_m}(x_*, \delta_*, \sigma_*)\,.
\end{equation}
Поскольку отображения~$f_m$ удовлетворяют соотношению~(\ref{eq2*A})
в $f_m(D),$ из~(\ref{eq6A}) получаем, что
\begin{equation}\label{eq14A}
M(\Gamma^{\,*}_m)\leqslant
(4N_0/l_0^n+(1/(\sigma_*-\delta_*))^n)\Vert Q\Vert_1:=c<\infty\,.
\end{equation}
Покажем, что соотношение~(\ref{eq14A}) противоречит условию слабой
плоскости отображённой области. В самом деле, по построению
$$h(|D^1_m|)\geqslant h(x_m, b_m^1) \geqslant
(1/2)\cdot h(f^{\,-1}_m(A), \partial D)>\delta/2\,,$$
\begin{equation}\label{eq14}
h(|D^2_m|)\geqslant h(x^{\,\prime}_m, b_m^2) \geqslant (1/2)\cdot
h(f^{\,-1}_m(A), \partial D)>\delta/2
\end{equation}
при всех $m\geqslant M_0$ и для некоторого $M_0\in {\Bbb N}.$
Положим $U:=B_h(x_0, r_0)=\{y\in \overline{{\Bbb R}^n}: h(y,
x_0)<r_0\},$ где $0<r_0<\delta/4$ и число $\delta$ относится к
соотношению~(\ref{eq14}). Заметим, что $|D^1_m|\cap
U\ne\varnothing\ne |D^1_m|\cap (D\setminus U)$ для каждого
$m\in{\Bbb N},$ поскольку $h(|D^1_m|)\geqslant \delta/2$ и $x_m\in
|D^1_m|,$ $x_m\rightarrow x_0$ при $m\rightarrow\infty.$ Аналогично,
$|D^2_m|\cap U\ne\varnothing\ne |D^2_m|\cap (D\setminus U).$
Поскольку $|D^1_m|$ и $|D^2_m|$ -- континуумы,
\begin{equation}\label{eq8A}
|D^1_m|\cap \partial U\ne\varnothing, \quad |D^2_m|\cap
\partial U\ne\varnothing\,,
\end{equation}
см., напр.,~\cite[теорема~1.I.5.46]{Ku}. Поскольку $\partial D$
слабо плоская, то для каждого $P>0$ найдётся окрестность $V\subset
U$ точки $x_0$ такая, что
\begin{equation}\label{eq9A}
M(\Gamma(E, F, D))>P
\end{equation}
для любых континуумов $E, F\subset D$ таких, что $E\cap
\partial U\ne\varnothing\ne E\cap \partial V$ и $F\cap \partial
U\ne\varnothing\ne F\cap \partial V.$ Покажем, что для достаточно
больших $m\in {\Bbb N}$
\begin{equation}\label{eq10A}
|D^1_m|\cap \partial V\ne\varnothing, \quad |D^2_m|\cap
\partial V\ne\varnothing\,.\end{equation}
В самом деле, $x_m\in |D^1_m|$ и $x^{\,\prime}_m\in |D^2_m|,$ где
$x_m, x^{\,\prime}_m\rightarrow x_0\in V$ при $m\rightarrow\infty.$
В таком случае, $|D^1_m|\cap V\ne\varnothing\ne |D^2_m|\cap V$ для
достаточно больших $m\in {\Bbb N}.$ Заметим, что $h(V)\leqslant
h(U)\leqslant 2r_0<\delta/2.$ В силу~(\ref{eq14})
$h(|D^1_m|)>\delta/2.$ Следовательно, $|D^1_m|\cap (D\setminus
V)\ne\varnothing$ и, значит, $|D^1_m|\cap\partial V\ne\varnothing$
(см., напр.,~\cite[теорема~1.I.5.46]{Ku}). Аналогично,
$h(V)\leqslant h(U)\leqslant 2r_0<\delta/2.$ Из~(\ref{eq14})
вытекает, что $h(|D^2_m|)>\delta,$ следовательно, $|D^2_m|\cap
(D\setminus V)\ne\varnothing.$ По~\cite[теорема~1.I.5.46]{Ku} мы
получим, что $|D^2_m|\cap\partial V\ne\varnothing.$ Таким образом,
соотношение~(\ref{eq10A}) установлено. Объединяя
соотношения~(\ref{eq8A}), (\ref{eq9A}) и (\ref{eq10A}), мы получим,
что $M(\Gamma^{\,*}_m)=M(\Gamma(|D^1_m|, |D^2_m|, D))>P.$ Последнее
условие противоречит~(\ref{eq14A}), что и доказывает теорему.~$\Box$

\medskip
{\bf 5. Примеры.}

\medskip
\begin{example}\label{ex1}
Рассмотрим в единичном шаре ${\Bbb B}^n$ последовательность
$f_m(z)=mz,$ $m=1,2,\ldots .$ Заметим, что как прямые отображения
$f_m,$ так и обратные отображения $f^{\,-1}_m$ удовлетворяют
условию~(\ref{eq2*A}) в соответствующей области при $Q\equiv 1$
(см., напр.,~\cite[теоремы~8.1, 8.5]{MRSY}). Заметим, что
последовательность $f_m$ не является равностепенно непрерывной, что
поясняется неинтегрируемостью функции $Q\equiv 1$ в ${\Bbb R}^n.$ В
то же время, обратная последовательность $f_m^{\,-1}:{\Bbb
R}^n\rightarrow {\Bbb R}^n$ равностепенно непрерывна в ${\Bbb R}^n,$
что непосредственно следует из теоремы~\ref{th1}. В самом деле,
рассмотрим произвольную точку $y_0\in{\Bbb R}^n$ и рассмотрим
сужение $g_m:=f^{\,-1}_m|_{B(y_0, r_0)},$ $r_0>0.$ Тогда при
некотором достаточно большом $m_0\in {\Bbb N}$ имеем: $g_m: B(y_0,
r_0)\rightarrow {\Bbb B}^n$ при $m\geqslant m_0.$ Если теперь
положить $Q(x)\equiv 1$ при $x\in {\Bbb B}^n$ и $Q(x)\equiv 0$ при
$x\not\in {\Bbb B}^n,$ то видно, что $g_m$
удовлетворяет~(\ref{eq2*A}) при $Q\in L^1({\Bbb R}^n).$ Таким
образом, все условия теоремы~\ref{th1} выполняются.
\end{example}

\medskip
\begin{example}\label{ex2}
В работах~\cite{SevSkv$_1$} и~\cite{SevSkv$_2$} были построены
примеры гомеоморфизмов, удовлетворяющих условию~(\ref{eq2*A}).
Укажем теперь на аналогичный пример семейства отображений с точками
ветвления. Пусть $p\geqslant 1$ настолько велико, что число $2/p$
меньше 1, и пусть, кроме того, $\alpha\in (0, 2/p)$ -- произвольное
число. Определим последовательность отображений $f_m: B(0,
2)\rightarrow {\Bbb B}^2$ следующим образом: $f_m(z)=(f\circ
\widetilde{f}_m)(z),$ $f(z)=z^2,$
$$\widetilde{f}_m(z)\,=\,\left
\{\begin{array}{rr} \frac{(|z|-1)^{1/\alpha}}{|z|}\cdot z\,, & 1+1/(m^{\alpha})\leqslant|z|\leqslant 2, \\
\frac{1/m}{1+(1/m)^{\alpha}}\cdot z\,, & 0<|z|<1+1/(m^{\alpha}) \ .
\end{array}\right.
$$
Используя подход, задействованный при
рассмотрении~\cite[предложение~6.3]{MRSY}, можно показать, что
внешняя дилатация $K_O(z, f_m)$ отображения $f_m$ в точке $z\in{\Bbb
B}^2$ вычисляется как $K_O(z,
f_m)=\frac{1}{\alpha}\cdot\frac{|z|}{|z|-1}$ при $z\in B(0,
2)\setminus \overline{{\Bbb B}^2},$ $K_O(z, f_m)=1$ при $z\in {\Bbb
B}^2.$ Заметим, что каждая точка $z\in {\Bbb B}^2\setminus B(0,
1/m)$ имеет ровно два прообраза при отображении $f_m,$ а именно,
$w=re^{i\varphi}\mapsto \pm((\sqrt{r})^{\,\alpha}+1)e^{i\varphi/2}.$
Точки $w\in B(0, 1/m)$ также имеют некоторые два прообраза при
отображении $f_m,$ скажем, $z^1_m$ и $z^2_m.$ Эти прообразы
совпадают и равны нулю только при $w=0.$
Таким образом, при $w\in {\Bbb B}^2$ и по выбору $\alpha$
$$K_I(w, f_m^{\,-1}, {\Bbb B}^2):=\sum\limits_{z\in f_m^{\,-1}(w)\cap {\Bbb B}^2}K_O(z, f_m)
\leqslant\frac{2}{\alpha}\cdot
\frac{(\sqrt{r})^{\alpha}+1}{(\sqrt{r})^{\alpha}}\leqslant\frac{4}{\alpha
r^{\alpha}}\in L^p({\Bbb B}^2)\,.$$

По~\cite[теоремы~8.1, 8.5]{MRSY} отображения $f_m$ удовлетворяют
соотношению~(\ref{eq2*A}) для функции
$Q(y)=\frac{4}{\alpha|y|^{\alpha}},$ $y\in {\Bbb B}^2;$ $Q(y)\equiv
0$ при $y\in {\Bbb R}^2\setminus \overline{{\Bbb B}^2}.$ Отсюда
следует, что для семейства отображения $\{f_m\}_{m=1}^{\infty}$
выполнены все условия теоремы~\ref{th1}. Также для этого семейства
выполнены все условия теоремы~\ref{th2}, поскольку $f_m$ фиксируют
бесконечное число точек единичного круга ${\Bbb B}^2$ при всех
$m\geqslant 2.$
\end{example}

\medskip
\begin{example}\label{ex3}
Теперь рассмотрим пример семейства отображений, аналогичный
примеру~\ref{ex2}, относящийся к пространству произвольной
размерности $n\geqslant 2.$ Для этого, для числа $p\geqslant 1$
такого, что $n/p(n-1),$ зафиксируем произвольны образом $\alpha\in
(0, n/p(n-1)).$ Пусть $f(x)=(r\cos 2\varphi, r\sin 2\varphi,
x_3,\ldots, x_n),$ $x=(x_1,x_2,\ldots, x_n),$
$r=\sqrt{x_1^2+x_2^2}.$ Определим последовательность отображений
$f_m: B(0, 2)\rightarrow {\Bbb B}^n$ следующим образом:
$f_m(x)=(f\circ \widetilde{f}_m)(x),$ где
$$\widetilde{f}_m(x)\,=\,\left
\{\begin{array}{rr} \frac{(|x|-1)^{1/\alpha}}{|x|}\cdot x\,, & 1+1/(m^{\alpha})\leqslant|x|\leqslant 2, \\
\frac{1/m}{1+(1/m)^{\alpha}}\cdot x\,, & 0<|x|<1+1/(m^{\alpha}) \ .
\end{array}\right.
$$
Используя подход, задействованный при
рассмотрении~\cite[предложение~6.3]{MRSY}, можно показать, что
$K_O(x,
f_m)\leqslant\frac{2^{n-1}}{\alpha^{n-1}}\cdot\frac{|x|^{n-1}}{(|x|-1)^{n-1}}$
при $x\in B(0, 2)\setminus \overline{{\Bbb B}^n},$ $K_O(x,
f)=2^{n-1}$ при $x\in {\Bbb B}^n.$ Заметим, что каждая точка $x\in
{\Bbb B}^n\setminus B(0, 1/m)$ имеет ровно два прообраза при
отображении $f_m,$ а именно,
$$w=(r\cos\varphi, r\sin\varphi,x_3,\ldots,x_n)\mapsto
\frac{1+|x|^{\alpha}}{|x|}\cdot (\pm r\cos\varphi/2, \pm
r\sin\varphi/2, x_3,\ldots, x_n)\,.$$
В таком случае,
$$K_I(w, f_m^{\,-1}, {\Bbb B}^n):=\sum\limits_{x\in f_m^{\,-1}(w)\cap {\Bbb B}^n}K_O(x, f_m)
\leqslant\frac{2^{n-1}}{\alpha^{n-1}}\cdot
\frac{(r^{\alpha}+1)^{n-1}}{(r^{\alpha(n-1)}}\leqslant
\frac{2^{n-1}4}{\alpha^{n-1} r^{\alpha(n-1)}}\in L^p({\Bbb
B}^n)\,.$$
По~\cite[теоремы~8.1, 8.5]{MRSY} отображения $f_m$ удовлетворяют
соотношению~(\ref{eq2*A}) для функции
$Q(y)=\frac{2^{n-1}4}{\alpha^{n-1}|y|^{(n-1)\alpha}},$ $y\in {\Bbb
B}^n;$ $Q(y)\equiv 0$ при $y\in {\Bbb R}^n\setminus \overline{{\Bbb
B}^n}.$ Отсюда следует, что для семейства отображения
$\{f_m\}_{m=1}^{\infty}$ выполнены все условия теоремы~\ref{th1}.
Также для этого семейства выполнены все условия теоремы~\ref{th2},
поскольку $f_m$ фиксируют бесконечное число точек единичного круга
${\Bbb B}^n$ при всех $m\geqslant 2.$
\end{example}

\medskip
{\bf 6. Лемма о континууме.} Вариант приводимого ниже утверждения
установлен в случае гомеоморфизмов в~\cite[пункт~v
леммы~2]{SevSkv$_1$}, см. также \cite[лемма~4.1]{SevSkv$_2$}. В
настоящей работе мы имеем дело с аналогичным утверждением,
относящимся к отображениям с ветвлением.

\medskip
\begin{lemma}\label{lem3}
{\sl\,Пусть $n\geqslant 2,$ $D$ и $D^{\,\prime}$ -- области в ${\Bbb
R}^n,$ при этом, $D$ имеет слабо плоскую границу, ни одна компонента
связности которой не вырождается в точку, а область $D^{\,\prime}$
локально связна на своей границе. Пусть также $A$ -- невырожденный
континуум в $D^{\,\prime}$ и $\delta>0.$ Предположим, $f_m$ --
последовательность открытых, дискретных и замкнутых отображений
области $D$ на $D^{\,\prime}$ со следующим свойством: для каждого
$m=1,2,\ldots$ найдётся континуум $A_m\subset D,$ $m=1,2,\ldots ,$
такой, что $f_m(A_m)=A$ и $h(A_m)\geqslant \delta>0.$ Если каждое
отображение $f_m$ удовлетворяет соотношению~(\ref{eq2*A}) в каждой
точке $y_0\in D^{\,\prime}$ и, при этом, $Q\in L^1(D^{\,\prime}),$
то найдётся $\delta_1>0$ такое, что
$$h(A_m,
\partial D)>\delta_1>0\quad \forall\,\, m\in {\Bbb
N}\,,$$
где $h(A_m,
\partial D)=\inf\limits_{x\in A_m, y\in \partial D}h(x, y).$}
\end{lemma}

\medskip
\begin{proof}
В силу компактности пространства~$\overline{{\Bbb R}^n}$ граница
области $D$ не пуста и является компактом, так что расстояние
$h(A_m,
\partial D)$ корректно определено.

\medskip
Проведём доказательство от противного. Предположим, что заключение
леммы не является верным. Тогда для каждого $k\in {\Bbb N}$ найдётся
номер $m=m_k$ такой, что $h(A_{m_k},
\partial D)<1/k.$ Можно считать, что последовательность $m_k$ возрастает по $k.$
Поскольку $A_{m_k}$ -- компакт, то найдутся $x_k\in A_{m_k}$ и
$y_k\in
\partial D$ такие, что $h(A_{m_k},
\partial D)=h(x_k, y_k)<1/k$ (см. рисунок~\ref{fig3}).
\begin{figure}[h]
\centerline{\includegraphics[scale=0.7]{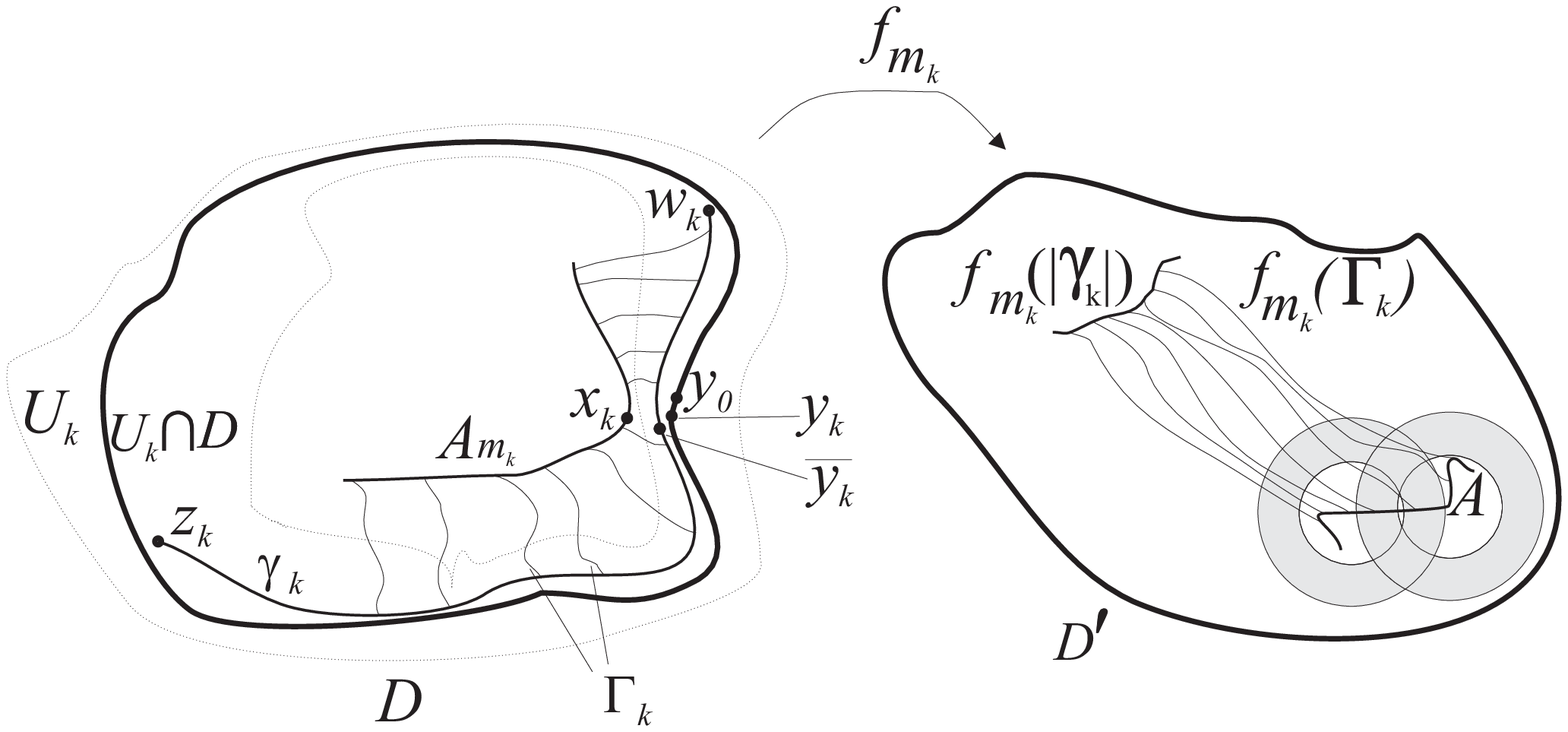}} \caption{К
доказательству леммы~\ref{lem3}}\label{fig3}
\end{figure}
Поскольку $\partial D$ -- компактное множество, мы можем считать,
что $y_k\rightarrow y_0\in
\partial D$ при $k\rightarrow \infty;$ тогда также
%
$x_k\rightarrow y_0\in \partial D$ при $k\rightarrow \infty.$
%
Пусть $K_0$ -- компонента связности $\partial D,$ содержащая точку
$y_0.$ Очевидно, $K_0$ -- континуум в $\overline{{\Bbb R}^n}.$
Поскольку $\partial D$ слабо плоская, а $D^{\,\prime}$ локально
связна на своей границе, то по теореме~\ref{th3} отображение
$f_{m_k}$ имеет непрерывное продолжение
$\overline{f}_{m_k}:\overline{D}\rightarrow
\overline{D^{\,\prime}}.$ Более того, отображение
$\overline{f}_{m_k}$ является равномерно непрерывным на
$\overline{D}$ при каждом фиксированном $k,$ так как
$\overline{f}_{m_k}$ непрерывно на компакте $\overline{D}.$ Тогда
для каждого $\varepsilon>0$ найдётся
$\delta_k=\delta_k(\varepsilon)<1/k$ такое, что
\begin{equation}\label{eq3B}
|\overline{f}_{m_k}(x)-\overline{f}_{m_k}(x_0)|<\varepsilon \quad
\forall\,\, x,x_0\in \overline{D},\quad h(x, x_0)<\delta_k\,, \quad
\delta_k<1/k\,.
\end{equation}
Выберем $\varepsilon>0$ таким, чтобы
\begin{equation}\label{eq5D}
\varepsilon<(1/2)\cdot {\rm  dist}\,(\partial D, A)\,.
\end{equation}
Обозначим $B_h(x_0, r)=\{x\in \overline{{\Bbb R}^n}: h(x, x_0)<r\}.$
Для фиксированного $k\in {\Bbb N},$ положим
$$B_k:=\bigcup\limits_{x_0\in K_0}B_h(x_0, \delta_k)\,,\quad k\in {\Bbb N}\,.$$
Поскольку $B_k$ -- окрестность континуума $K_0,$
по~\cite[лемма~2.2]{HK} найдётся окрестность $U_k$ множества $K_0,$
такая, что $U_k\subset B_k$ и $U_k\cap D$ связно. Можно считать, что
$U_k$ -- открыто, так что $U_k\cap D$ является линейно связным
(см.~\cite[предложение~13.1]{MRSY}). Пусть $h(K_0):=\sup\limits_{x,
y\in K_0}h(x, y)=m_0.$ Тогда найдутся $z_0, w_0\in K_0$ такие, что
$h(K_0)=h(z_0, w_0)=m_0.$ Следовательно, найдутся последовательности
$\overline{y_k}\in U_k\cap D,$ $z_k\in U_k\cap D$ и $w_k\in U_k\cap
D$ такие, что $z_k\rightarrow z_0,$ $\overline{y_k}\rightarrow y_0$
и $w_k\rightarrow w_0$ при $k\rightarrow\infty.$ Можно считать, что
\begin{equation}\label{eq2B}
h(z_k, w_k)>m_0/2\quad \forall\,\, k\in {\Bbb N}\,.
\end{equation}
Поскольку множество $U_k\cap D$ линейно связно, мы можем соединить
точки $z_k,$ $\overline{y_k}$ и $w_k,$ используя некоторую кривую
$\gamma_k\in U_k\cap D.$ Как обычно, мы обозначаем через
$|\gamma_k|$ носитель (образ) кривой $\gamma_k$ в области $D.$ Тогда
$f_{m_k}(|\gamma_k|)$ -- компактное множество в $D.$ Если
$x\in|\gamma_k|,$ то найдётся $x_0\in K_0$ такое, что $x\in B(x_0,
\delta_k).$ Положим $\omega\in A\subset D.$ Поскольку
$x\in|\gamma_k|$ и, более того, $x$ -- внутренняя точка $D,$ мы
можем использовать запись $f_{m_k}(x)$ вместо
$\overline{f}_{m_k}(x).$ Из соотношений~(\ref{eq3B}) и~(\ref{eq5D}),
а также по неравенству треугольника, мы получаем, что
$$|f_{m_k}(x)-\omega|\geqslant
|\omega-\overline{f}_{m_k}(x_0)|-|\overline{f}_{m_k}(x_0)-f_{m_k}(x)|\geqslant$$
\begin{equation}\label{eq4C}
\geqslant {\rm  dist}\,(\partial D^{\,\prime}, A)-(1/2)\cdot{\rm
dist}\,(\partial D^{\,\prime}, A)=(1/2)\cdot{\rm  dist}\,(\partial
D^{\,\prime}, A)>\varepsilon
\end{equation}
для достаточно больших $k\in {\Bbb N},$ где ${\rm  dist}\,(\partial
D, A):=\inf\limits_{x\in \partial D, y\in A}|x-y|.$ Переходя к
$\inf$ в~(\ref{eq4C}) по всем $x\in |\gamma_k|$ и $\omega\in A,$ мы
получаем, что
\begin{equation}\label{eq6B}
{\rm dist}\,(f_{m_k}(|\gamma_k|), A):=\inf\limits_{x\in
f_{m_k}(|\gamma_k|), y\in A}|x-y|>\varepsilon, \quad\forall\,\,
k=1,2,\ldots \,.
\end{equation}
Покроем множество $A$ шарами $B(x, \varepsilon/4),$ $x\in A.$
Поскольку $A$ компакт, мы можем считать, что $A\subset
\bigcup\limits_{i=1}^{M_0}B(x_i, \varepsilon/4),$ $x_i\in A,$
$i=1,2,\ldots, M_0,$ $1\leqslant M_0<\infty.$ По определению, $M_0$
зависит только от $A,$ в частности, $M_0$ не зависит от $k.$ Положим
\begin{equation}\label{eq5C}
\Gamma_k:=\Gamma(A_{m_k}, |\gamma_k|, D)\,.
\end{equation}
Пусть $\Gamma_{ki}:=\Gamma_{f_k}(x_i, \varepsilon/4,
\varepsilon/2),$ другими словами, $\Gamma_{ki}$ состоит из всех
кривых $\gamma:[0, 1]\rightarrow D,$ принадлежащих семейству
$\Gamma_k,$ таких что $\gamma(0)\in S(x_i, \varepsilon/4),$
$\gamma(1)\in S(x_i, \varepsilon/2)$ и $\gamma(t)\in A(x_i,
\varepsilon/4, \varepsilon/2)$ при $0<t<1.$ Покажем, что
\begin{equation}\label{eq6C}
\Gamma_k>\bigcup\limits_{i=1}^{M_0}\Gamma_{ki}\,.
\end{equation}
В самом деле, пусть $\widetilde{\gamma}\in \Gamma_k.$ Тогда
$\gamma:=f_{m_k}(\gamma)\in \Gamma(A, f_{m_k}(|\gamma_k|),
D^{\,\prime}).$ Поскольку шары $B(x_i, \varepsilon/4),$ $1\leqslant
i\leqslant M_0,$ образуют покрытие компакта $A,$ найдётся $i\in
{\Bbb N}$ такое, что $\gamma:[0, 1]\rightarrow D,$ $\gamma(0)\in
B(x_i, \varepsilon/4)$ и $\gamma(1)\in f_{m_k}(|\gamma_k|).$ По
соотношению~(\ref{eq6B}), $|\gamma|\cap B(x_i,
\varepsilon/4)\ne\varnothing\ne |\gamma|\cap (D^{\,\prime}\setminus
B(x_i, \varepsilon/4)).$ Следовательно,
по~\cite[теорема~1.I.5.46]{Ku} найдётся $0<t_1<1$ такое, что
$\gamma(t_1)\in S(x_i, \varepsilon/4).$ Можно считать, что
$\gamma(t)\not\in B(x_i, \varepsilon/4)$ при $t>t_1.$ Положим
$\gamma_1:=\gamma|_{[t_1, 1]}.$ Из~(\ref{eq6B}) вытекает, что
$|\gamma_1|\cap B(x_i, \varepsilon/2)\ne\varnothing\ne
|\gamma_1|\cap (D\setminus B(x_i, \varepsilon/2)).$ Следовательно,
по~\cite[теорема~1.I.5.46]{Ku} найдётся $t_1<t_2<1$ такое, что
$\gamma(t_2)\in S(x_i, \varepsilon/2).$ Можно считать, что
$\gamma(t)\in B(x_i, \varepsilon/2)$ при всех $t<t_2.$ Полагая
$\gamma_2:=\gamma|_{[t_1, t_2]},$ заметим, что кривая $\gamma_2$
является подкривой $\gamma,$ принадлежащей семейству $\Gamma(S(x_i,
\varepsilon/4), S(x_i, \varepsilon/2), A(x_i, \varepsilon/4,
\varepsilon/2)).$ Таким образом, соотношение~(\ref{eq6C})
установлено. Положим
$$\eta(t)= \left\{
\begin{array}{rr}
4/\varepsilon, & t\in [\varepsilon/4, \varepsilon/2],\\
0,  &  t\not\in [\varepsilon/4, \varepsilon/2]\,.
\end{array}
\right. $$
Заметим, что $\eta$ удовлетворяет соотношению~(\ref{eqA2}) при
$r_1=\varepsilon/4$ и $r_2=\varepsilon/2.$ Поскольку отображение $f$
удовлетворяет соотношению~(\ref{eq2*A}), то полагая здесь $y_0=x_i,$
получаем:
\begin{equation}\label{eq8C}
M(\Gamma_{f_k}(x_i, \varepsilon/4, \varepsilon/2))\leqslant
(4/\varepsilon)^n\cdot\Vert Q\Vert_1<c<\infty\,,
\end{equation}
where $c$ -- некоторая положительная постоянная и $\Vert Q\Vert_1$
--  $L_1$-норма функции $Q$ в $D^{\,\prime}.$ Из~(\ref{eq6C}) и
(\ref{eq8C}), учитывая полуаддитивность модуля семейств кривых,
получаем:
\begin{equation}\label{eq4B}
M(\Gamma_k)\leqslant
\frac{4^nM_0}{\varepsilon^n}\int\limits_{D^{\,\prime}}Q(y)\,dm(y)\leqslant
c\cdot M_0<\infty\,.
\end{equation}
С другой стороны, рассуждая также, как при доказательстве
соотношения~(\ref{eq7}) и используя условие~(\ref{eq2B}), мы
получим, что $M(\Gamma_k)\rightarrow\infty$ при
$k\rightarrow\infty,$ то противоречит~(\ref{eq4B}). Полученное
противоречие доказывает лемму.~$\Box$
\end{proof}

\medskip
{\bf 7. Об отображениях с фиксированной точкой.} Отметим, что в
теореме~\ref{th2} присутствуют достаточно жёсткие условия как на
границу отображённой области, так и на само семейство. Сейчас мы
укажем некоторый частный случай, в котором указанные условия могут
быть сформулированы в более изящном виде. Для этой цели рассмотрим
следующее определение. Для областей $D, D^{\,\prime}\subset {\Bbb
R}^n,$ точек $a\in D,$ $b\in D^{\,\prime}$ и произвольной измеримой
по Лебегу функции $Q:D^{\,\prime}\rightarrow [0, \infty]$ обозначим
через ${\frak S}_{a, b, Q}(D, D^{\,\prime})$ семейство всех
открытых, дискретных и замкнутых отображений $f$ области $D$ на
$D^{\,\prime},$ удовлетворяющих условию~(\ref{eq2*A}) для каждого
$y_0\in f(D),$ таких что $f(a)=b.$ Справедливо следующее
утверждение.

\medskip
\begin{theorem}\label{th4}
{\sl Предположим, что область $D$ имеет слабо плоскую границу, ни
одна из связных компонент которой невырождена. Если $Q\in
L^1(D^{\,\prime})$ и область $D^{\,\prime}$ локально связна на своей
границе, то $f\in {\frak S}_{a, b, Q}(D, D^{\,\prime})$ продолжается
по непрерывности до отображения
$\overline{f}:\overline{D}\rightarrow \overline{D^{\,\prime}},$
$\overline{f}|_D=f,$ при этом,
$\overline{f}(\overline{D})=\overline{D^{\,\prime}}$ и семейство
${\frak S}_{a, b, Q}(\overline{D}, \overline{D^{\,\prime}}),$
состоящее из всех продолженных отображений
$\overline{f}:\overline{D}\rightarrow \overline{D^{\,\prime}},$
равностепенно непрерывно в $\overline{D}.$ }
\end{theorem}

\medskip
\begin{proof} Равностепенная непрерывность семейства ${\frak S}_{a, b, Q}(D,
D^{\,\prime}),$ возможность непрерывного продолжения на границу
каждого $f\in{\frak S}_{a, b, Q}(D, D^{\,\prime})$ и равенство
$\overline{f}(\overline{D})=\overline{D^{\,\prime}}$ вытекают из
теорем~\ref{th1} и~\ref{th3}. Осталось установить равностепенную
непрерывность семейства продолженных отображений
$\overline{f}:\overline{D}\rightarrow \overline{D^{\,\prime}}$ в
точках границы области~$D.$

\medskip
Докажем это утверждение от противного. Предположим, что семейство
${\frak S}_{a, b, Q}(\overline{D}, \overline{D^{\,\prime}})$ не
является равностепенно непрерывным в некоторой точке $x_0\in\partial
D.$ Тогда найдутся точки $x_m\in D$ и отображения $f_m\in {\frak
S}_{a, b, Q}(\overline{D}, \overline{D^{\,\prime}}),$ $m=1,2,\ldots
,$ такие что $x_m\rightarrow x_0$ при $m\rightarrow\infty$ и, при
этом, при некотором $\varepsilon_0>0$
\begin{equation}\label{eq15}
h(f_m(x_m), f_m(x_0))\geqslant\varepsilon_0\,,\quad m=1,2,\ldots\,.
\end{equation}
Выберем произвольным образом точку $y_0\in D^{\,\prime},$ $y_0\ne
b,$ и соединим её с точкой $b$ некоторой кривой в $D^{\,\prime},$
которую мы обозначим $\alpha.$ Положим $A:=|\alpha|.$ Пусть $A_m$ --
полное поднятие кривой $\alpha$ при отображении $f_m$ с началом в
точке $a$ (оно существует ввиду~\cite[лемма~3.7]{Vu$_1$}). Заметим,
что $h(A_m, \partial D)>0$ по замкнутости отображения~$f_m.$ Теперь
возможны следующие случаи: либо $h(A_m):=\sup\limits_{x, y\in
A_m}h(x, y)\rightarrow 0$ при $m\rightarrow\infty,$ либо
$h(A_{m_k})\geqslant\delta_0>0$ при $k\rightarrow\infty$ для
некоторой возрастающей последовательности номеров $m_k$ и некоторого
$\delta_0>0.$

В первом из этих случаев, очевидно, $h(A_m, \partial D)\geqslant
\delta>0$ при некотором $\delta>0.$ Тогда семейство отображений
$\{f_m\}_{m=1}^{\infty}$ равностепенно непрерывно в точке $x_0$ по
теореме~\ref{th2}, что противоречит условию~(\ref{eq15}).

Во втором случае, если $h(A_{m_k})\geqslant\delta_0>0$ при
$k\rightarrow\infty,$ мы также имеем, что $h(A_{m_k}, \partial
D)\geqslant \delta_1>0$ при некотором $\delta_1>0$ в силу
леммы~\ref{lem3}. Опять же, по теореме~\ref{th2} мы имеем, что
семейство $\{f_{m_k}\}_{k=1}^{\infty}$  равностепенно непрерывно в
точке $x_0,$ и это противоречит условию~(\ref{eq15}).

Итак, в обоих из двух возможных случаев мы пришли к противоречию
с~(\ref{eq15}), и это указывает на неверность предположения об
отсутствии равностепенной непрерывности семейства ${\frak S}_{a, b,
Q}(D, D^{\,\prime})$ в $\overline{D}.$ Теорема доказана.~$\Box$
\end{proof}

\medskip
{\bf Открытый вопрос.} В условиях теоремы~\ref{th4} участвует
условие невырожденности произвольной связной компоненты границы
области $D,$ и оно существенно используется при доказательстве. {\bf
Верно ли аналогичное заключение без этого условия ?}


КОНТАКТНАЯ ИНФОРМАЦИЯ

\medskip
\noindent{{\bf Евгений Александрович Севостьянов} \\
{\bf 1.} Житомирский государственный университет им.\ И.~Франко\\
кафедра математического анализа, ул. Большая Бердичевская, 40 \\
г.~Житомир, Украина, 10 008 \\
{\bf 2.} Институт прикладной математики и механики
НАН Украины, \\
отдел теории функций, ул.~Добровольского, 1 \\
г.~Славянск, Украина, 84 100\\
e-mail: esevostyanov2009@gmail.com}

\medskip
\noindent{{\bf Сергей Александрович Скворцов, Александр Петрович Довгопятый} \\
Житомирский государственный университет им.\ И.~Франко\\
кафедра математического анализа, ул. Большая Бердичевская, 40 \\
г.~Житомир, Украина, 10 008 \\ e-mail: serezha.skv@gmail.com }

\end{document}